# Emergence of dynamic properties in network hyper-motifs

Miri Adler[1] and Ruslan Medzhitov[2,*]

[1]Broad Institute of Massachusetts Institute of Technology and Harvard, Cambridge, MA, USA

[2]Howard Hughes Medical Institute, Department of Immunobiology, Yale University School of Medicine, New Haven, CT 06510, USA

*Corresponding author. Email: ruslan.medzhitov@yale.edu

**Abstract:**

Networks are fundamental for our understanding of complex systems. Interactions between individual nodes in networks generate network motifs - small recurrent patterns that can be considered the network's building-block components, providing certain dynamical properties. However, it remains unclear how network motifs are arranged within networks and what properties emerge from interactions between network motifs. Here we develop a framework to explore the mesoscale-level behavior of complex networks. Considering network motifs as hypernodes, we define the rules for their interaction at the network's next level of organization. We infer the favorable arrangements of interactions between network motifs into hyper-motifs from real evolved and designed networks data including biological, neuronal, social, linguistic and electronic networks. We mathematically explore the emergent properties of these higher-order circuits and their relations to the properties of the individual minimal circuit components they combine. This framework provides a basis for exploring the mesoscale structure and behavior of complex systems where it can be used to reveal intermediate patterns in complex networks and to identify specific nodes and links in the network that are the key drivers of the network's emergent properties.

## Introduction

Complex systems describe a collection of multiple agents that influence each other according to some specified rules. A common and powerful tool in exploring the structure of complex systems is the use of a network description, where the nodes represent the individual agents and the edges represent the interactions between them (*1*, *2*). The study of networks has uncovered common principles that underlie the behavior of vastly different fields of study, including physics, biology, sociology and engineering (*3–6*). One of these common principles is the existence of network motifs (*7*). Network motifs are small recurrent patterns that can provide certain features that are important for the specific network. Network motifs such as the feedforward loop (FFL) have been studied separately both theoretically and experimentally in various fields and their dynamical properties were elucidated (*8–10*). However, it remains unclear how network motifs are joined in real networks to make larger circuits and what properties can emerge from these higher-level functional modules (Fig. 1A).

Indeed, it has been established that although the same network motifs appear in different contexts, the way that they are joined varies to provide distinct features. For example, FFLs are joined in bacterial transcription networks with multiple outputs to produce a first-in-first-out (FIFO) response. During sporulation, *B. subtilis* uses cascades of FFLs to activate genes in a series of temporal waves. In neuronal networks, FFLs can combine with multiple inputs to provide coincidence detection and pain relief when multiple pain sensation inputs are integrated into a single output (*11*, *12*) (Fig. 1B). Thus, combining the same kind of canonical motifs in different ways can generate novel properties and can also be used to 'silence' certain motif properties when they are not needed.

Work on modularity (*13–15*) and networks of networks (*16*) has revealed hierarchical structures (*17*) where several levels of organization are sometimes needed to describe the network. Recently, Battiston et al. (*18*) reviewed existing approaches to explore higher-order interactions in complex networks including the use of hypergraphs and simplicial complexes (*19*), and highlighted the challenges in the field where higher-order interactions are hard to infer from real data that are mostly based on simple pairwise interactions. Previous studies have explored higher-order clusters and generalizations of network motifs (*20–23*), and studied network motifs at different scales in real networks (*24*). However, there are currently no approaches to explore how building-block circuits such as network motifs interact with each other in the network based on specific circuit topologies to form the next level of organization. In order to understand the origin of emergent properties in complex systems it is important to study how network motifs interact with each other as defining the rules of atom interactions that form molecules was crucial for the understanding of the macroscale behavior of matter.

Here, we develop a framework for exploring the rules for emergent properties at intermediate levels of organization of complex networks. In this framework, we consider that the individual agents that interact with each other are the minimal building-block circuit topologies in the network and explore the emergent properties that result from the way that they are embedded in the network, which we call hyper-motifs. We develop a method to explore the favorable arrangements of these hyper-motifs in real networks and apply it to biological, neuronal, social, linguistic and electronic networks. This approach sheds light on the inner structure of complex networks and reveals new levels of organizations in real evolved and designed networks.

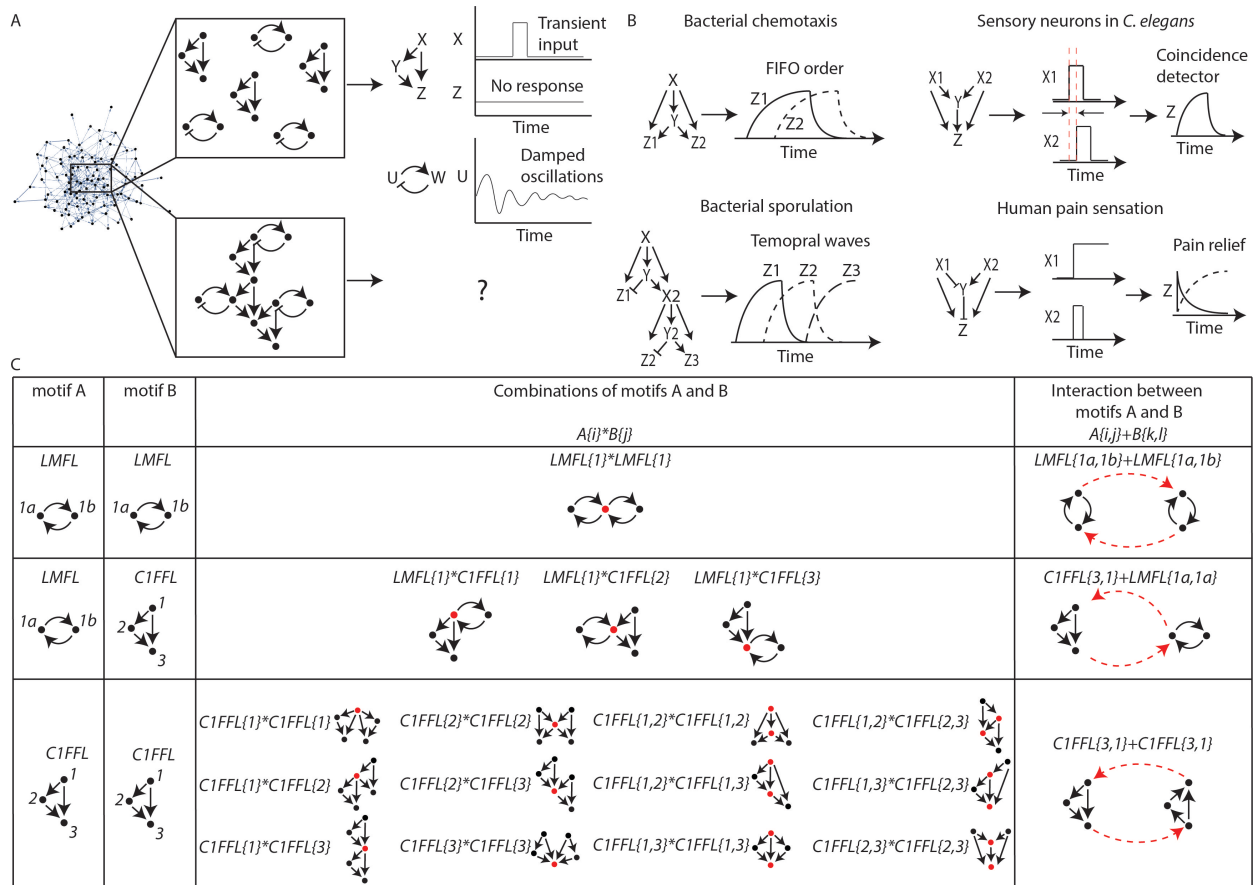

**Figure 1: A framework to explore hyper-motifs in complex networks. A.** Network motifs such as the coherent type 1 feedforward loop (C1FFL) and the mutual feedback loop where u is an activator and w is a repressor can filter out transient input signals and generate damped oscillations, respectively. But it is unclear what would be the properties of the circuits made up of combinations of different network motifs. **B.** Examples of four real networks in which the FFL is a network motif, but the way it is joined with other FFLs is different and provides distinct dynamical properties. **C.** A table that exemplifies the two ways in which network motifs can be directly joined together for the Lock-ON mutual feedback loop (LMFL) and the C1FFL. Combination of motifs A and B is where the two motifs share at least one vertex. These are marked by A{i}*B{j}, where {i}, {j} is the set of nodes each motif is sharing. Interaction of motifs A and B is where the two motifs are linked through at least one edge, A{i,j}+B{k,l}, where i and k are the set of sender nodes, and j and l are the set of receiver nodes from motifs A and B, respectively. The shared vertices in the 'Combinations' column and the linking edges in the 'Interactions' column are marked in red.

## Results

### A framework for exploring high-level modules of network motifs

To explore how two network motifs are integrated to form hyper-motifs - higher-level network modules, we consider two ways in which motifs can be directly joined. The first is where the two motifs share at least one node, which we define as a *combination* of motifs. The maximal number of shared nodes ($N_{V,max}$) must be smaller than the size of the smallest motif embedded such that the autonomy of each motif's topology is maintained. Therefore, $N_{V,max} = min(n_A, n_B) - 1$ where $n_A$ and $n_B$ are the sizes of motifs A and B, respectively. The identity of the nodes that are shared between the motifs is based on categorizing the nodes in each motif according to their unique roles. For example, nodes that participate in FFLs will be categorized into three distinct groups: input, intermediate and output nodes. If the roles in a certain motif are symmetric, we consider them in the same category (*21*). We mark a combination of motifs A and B in which nodes {*i*} of motif A are shared with nodes {*j*} of motif B as A{*i*}*B{*j*} (Fig. 1C).

The second way two network motifs can be joined, which we define as an *interaction* between motifs, is when they are linked through at least one edge. When two motifs interact, every pair of vertices that do not participate in the same motif can be linked. Therefore, the maximal number of linking edges is $N_{E,max} = 2n_A n_B$ for directed networks, and half as much for undirected networks. We mark an interaction between motifs A and B as A{*i,j*}+B{*k,l*}, where there are links from nodes {*i*} of motif A to nodes {*l*} of motif B and links from nodes {*k*} of motif B to nodes {*j*} of motif A (Fig. 1C).

In Figure 1C we exemplify these definitions for several pairs of motifs. Two mutual feedback circuits can only be combined by sharing one vertex. A mutual feedback and an FFL can be combined by sharing one vertex which can be either the FFL's input, intermediate or output node. The space of possible combinations increases substantially as the size of the motifs that are being combined increases. For example, two FFLs have 12 different ways to be joined when they share either one or two vertices. Note that we list in Figure 1C the core topology of each possible combination. However, each such combination can be extended where every pair of vertices that do not participate in the same motif can be linked (Methods, fig. S1). For each pair of motifs in Figure 1C, we provide an example of the way they can interact (see Methods for all possible topologies of motif interactions, fig. S1). This framework thus provides a way to count all possible circuit topologies for combinations and interactions of two motifs.

**Observed combinations of network motifs in real networks**

To explore whether network motifs are joined in real networks in specific ways, we developed a method to detect enriched combinations of network motifs in large networks. Given a real network of interest, the first step is to detect the network motifs (up to size n) that characterize the network. There are several algorithms for finding network motifs in large networks (*7, 25, 26*). After the network motifs of up to n-order have been identified, we categorize the nodes that participate in the network motifs based on their role in the motifs. Next, we compute the level of overlap between nodes in every pair of motif roles. A large overlap between two groups of motif roles means that the motifs are often combined in the network by sharing these nodes (Fig. 2A). Finally, in order to test the statistical significance of these combinations, we compare the observed overlap in the real network to the overlap in randomized networks when we keep properties of the network including the degree distribution and the frequency of all subgraphs up to size n the same as in the real network. This comparison allows us to find statistically significant over- and under-represented combinations of network motifs that do not emerge due to topological constraints in the network (*27, 28*) (Methods).

We applied this method for several natural and designed networks of different origin. In the *E. coli* transcription network (*29*) the nodes are transcription factors and their target genes, and the edges represent regulatory interactions. In this network the self-loop and FFL motifs are often combined such that the intermediate node of the FFL shows autoregulation. In the neuronal network of *C. elegans* (*7, 30*) where edges represent synaptic connection between neurons, there are 6 network motifs including the FFL and 5 different versions of mutual feedback circuits. We find a large number of over-represented combinations of these motifs in the network, in line with the abundance of evidence for high-order circuits in neuronal systems (*31, 32*). The over-represented combinations show that the neuronal network has a layered structure where in most cases an output node of one motif serves as an input of another motif. Moreover, most combinations in which two motifs are not combined in a layered way are excluded in the network. One exception to this structure is a combination where two double mutual feedback circuits are combined not in a layered manner. We discuss potential emergent properties of this combined circuit in the next section. Interestingly, there is a large number of combinations that are excluded in the *C. elegans* neuronal network (fig. S2), supporting the notion that network motifs are not distributed randomly in the network but are arranged in a way that provides a desired functionality in a given system. In an electronic circuit network (digital fractional multipliers) (*33*, *34*) where the nodes represent different logic gates and flip-flops, the 3-node feedback loop circuits are often combined by sharing two nodes with other 3-node feedback loop circuits. In a food web of lizards on the St. Martin island (*35*) where the edges represent predator-prey relations, the network motifs are a 3-node chain and an FFL. These motifs are often combined in the network where they either share the highest predator or the lowest prey in the motifs, or both. The network motifs in this food web tend not to share intermediate-level predators or preys and avoid very long food chains. In a citation network (*36*, *37*) where nodes are researchers and edges represent scientific citations, the FFL motifs are arranged in cascades. This structure is in line with the fact that citation networks are a type of an information network where the flow in the network is possible only in one direction. The over- and under-represented combinations in the citation network further show that scientists tend to cite recent papers more than original ones. We considered a social network (*38*) where nodes are people and an edge is drawn from person A to person B if person A considers person B as a close friend. In this friendship network, the network motifs include the FFL and 4 versions of mutual feedback circuits. These motifs are often combined such that two FFLs share the input and intermediate nodes, and two regulating mutual feedback loops share the nodes that mutually interact with each other. In these enriched combinations the 'output' nodes usually show mutual links in the network. This suggests a pattern where popular people are liked by either members of the same clique or individuals that don't show a mutual friendship. Finally, we analyzed word adjacency networks of texts in English and Japanese where each node represents a word and a directed connection occurs when one word directly follows the other in the text (*39*). The network motifs in both languages are a 3-node cascade, a regulating and regulated V circuits. We find that these motifs show two distinct patterns where they either combine in a layered manner, or that certain words have multiple words that are adjacent to them (Fig. 2B, Methods). Interestingly, the examination of combinations of network motifs in these linguistic networks shows that there are common principles in the way the network motifs are integrated where the enriched and excluded combinations of network motifs in the Japanese text are also found in the English network. See figure S2 for information on the statistical significance of all over-represented combinations we detected in real networks.

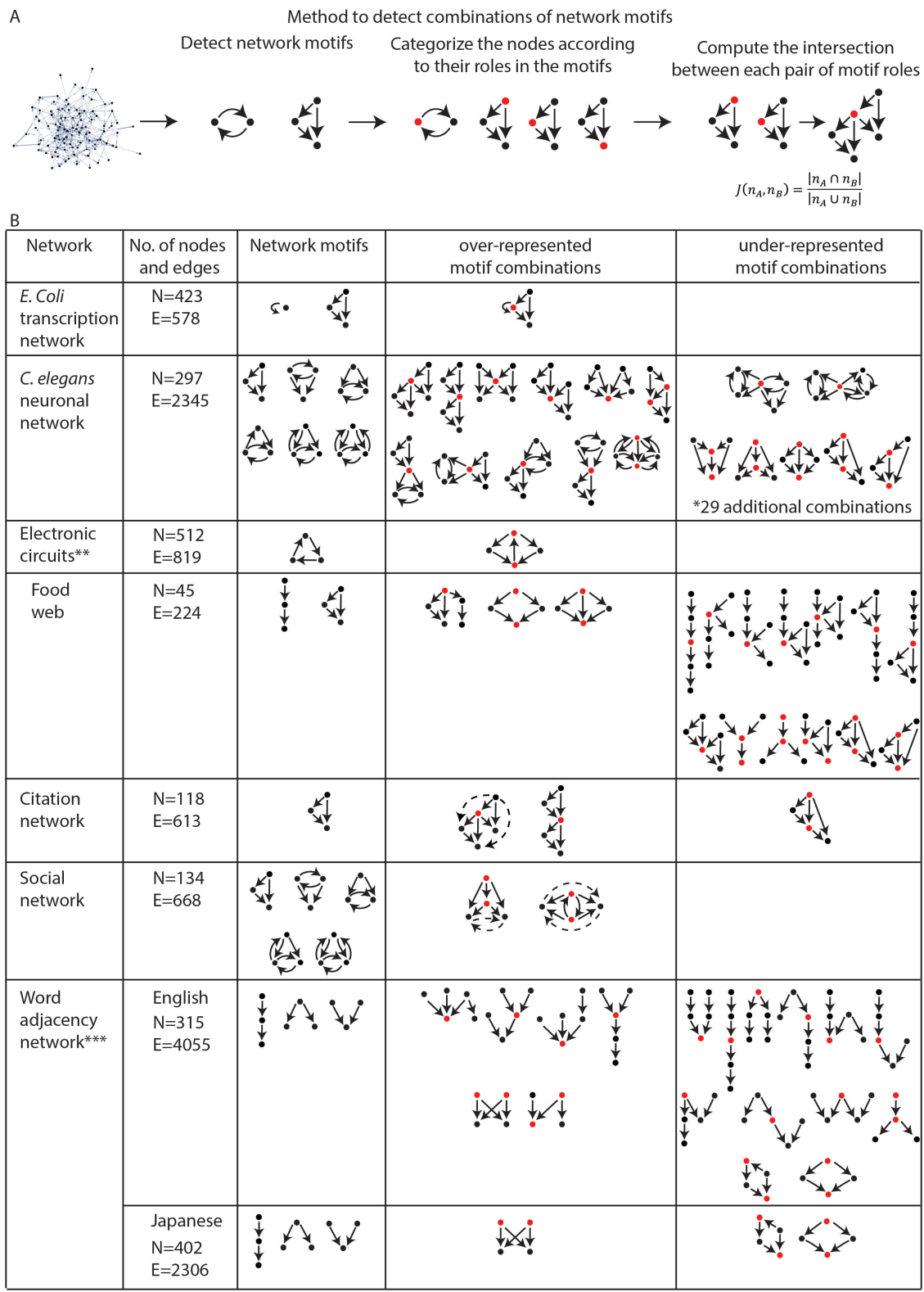

**Figure 2: A methods to detect enriched combinations of network motifs in real networks. A.** Schematics of the method to detect enriched combinations of network motifs in real networks. **B.** A table that summarizes the analysis of real networks with the type of network, its number of nodes and edges, its network motifs, and over- and under-represented combinations. Shared nodes are marked in red. We list in the SI the full list of under-represented combinations of the *C. elegans* neuronal network (fig. S2B). In the citation and social networks, the dashed arrows in the over-represented combinations represent additional edges that often appear in the network in addition to the core topology of the combined circuit (see fig. S2A for their frequencies in the networks). **We analyzed three different electronic circuits. The combined circuit in which two 3-node loop circuits share an edge is significantly enriched in all of them. The other two electronic circuits have N=122, E=189, and N=252, E=399. ***We downsampled the word adjacency networks (see Methods for details about the downsampling method we used).

**Emergent properties of combinations of building-block circuits**

We next explore the dynamical properties of combinations of network motifs. In our modeling framework, we use nonlinear Hill functions to describe relations between the nodes in a circuit (Methods). As our minimal building-block circuits, we consider three classes of circuit topologies that were found as network motifs in real natural and engineered networks (*11*). The first is a self-loop (SL) circuit which is a simple motif with only one node (X) positively autoregulating its own levels that can provide bistability. This means that X converges to a high steady state level only if its initial level is above a certain threshold, otherwise it declines to zero (Fig. 3A). The second class of circuits we consider includes three types of mutual feedback circuits: the toggle switch circuit (TMFL) in which X and Y mutually inhibit each other leading to a switch between their final levels, the Lock-ON circuit (LMFL) in which X and Y are both turned ON or OFF due to their mutual activation, and the oscillator circuit (OMFL) with X as a repressor and Y as an activator (Fig. 3A). The third class of circuits is the FFL circuit with two main types - a coherent type 1 FFL (C1FFL) in which the input (X) activates an intermediate node (Y) and both X and Y activate the output (Z), and an incoherent type 1 FFL (I1FFL) where Y is a repressor of Z. We show in Fig. 3A examples of previously explored dynamical features that the FFL circuits can exhibit for a wide range of parameters (Methods). We note that in general when we discuss a circuit's dynamical properties, we are referring to properties that are observed for a given choice of models (e.g., nonlinear relations, AND vs OR logic gates) and a certain range of model parameters. However, the circuits may show other properties for other choices of models or parameter values (*40*, *41*). When we discuss the properties of a combination of circuits, we compare them to the properties of each circuit component when we keep the same choice of parameter values.

Combining a self-loop motif with other motifs generally provides bistability to the node that is positively self-regulating. Therefore, a toggle-switch circuit combined with self-loops gives rise to a new stable state in which both X and Y are turned OFF (Fig. 3B). The oscillator circuit with both nodes self-regulating their own levels shows two additional stable steady states on top of the oscillating state - an OFF state where both X and Y are turned OFF, and a state where only Y (the activator) is turned ON. A state where only X (the repressor) is turned ON is not possible since the repressor cannot increase its levels without the presence of the activator Y (Fig. 3B). A Lock-ON circuit combined with self-loops does not provide new steady states since bistability is already provided in a simple Lock-ON circuit. Here we focus on a positive self-loop motif that can provide new steady states. A negative self-loop motif combined with other motifs does not provide bistability but can accelerate the circuits' response time (*42*).

When exploring combinations of a self-loop motif and FFL circuits, we asked whether the FFL's behavior is sensitive to the identity of the node that has a self-loop. To that aim, we compared the dynamics of the output Z of a simple FFL (without self-loops) to FFLs with a self-loop on the FFL's input, intermediate or output nodes. The FFL combined with self-loops provides bistability to the FFL's response where the final level of the output depends on the initial level of the input. However, FFL circuits in which the output node is autoregulated are less sensitive to the initial levels of the input (Fig. 3C). The incoherent FFL shows an interesting nontrivial behavior when the self-loop is on the intermediate repressor, Y. For high initial levels of input, the output shows a pulsatile behavior, which is similar to the behavior of a simple incoherent FFL. However, if the initial levels of input are low, the output rises to a high steady state level with a delayed response. The reason for this behavior is that for low levels of input, the repressor Y declines to zero due to its self-loop, thus lifting the repression on the output Z and allowing it to increase its levels (Fig. 3C). The emergence of this new high steady state level is possible only if the intermediate node is the one that is positively autoregulated (Methods). Thus, the self-loop provides an additional thresholding mechanism for the behavior of the dynamical circuits it is coupled with. This mathematical model therefore provides a possible explanation for why the FFL and self-loop motifs are combined in the *E. coli* transcription network such that the FFL's intermediate node has a self-loop interaction. This combination may provide bistability where the *E. coli* target genes are sensitive to the levels of the input signal.

We next show that different combinations of the same two motifs can yield different dynamical properties when keeping the same parameter values. To demonstrate this, we consider three different combinations of two oscillator circuits. Exploring the dynamical behavior of the combinations of circuits for varying initial conditions shows that the different combinations behave dynamically differently and converge to different steady states (Fig. 3D, Methods).

We find that coupling the oscillator circuit with the coherent and incoherent FFLs through their intermediate node yields interesting features of a pulsatile response (for a coherent FFL) or a delayed rising response (for an incoherent FFL) for low and high initial levels of W. The width of the pulse and the duration of the delay are proportional to the initial levels of W (Fig. 3E, Methods). In the supplementary information we show the modeling results of the combinations of all possible pairs of circuits from these three classes of circuit topologies (fig. S3, S4).

There are several observations that this framework demonstrates when one considers the potential emergent properties of combinations of network motifs. First, the way that the motifs are combined or the identity of the shared nodes that link the two motifs is an important factor that can drastically affect the resulting behavior. Second, we find that when multistablity emerges from combining two motifs, oftentimes the combined circuit can both preserve the autonomy of each circuit component where it shows their individual properties, and show emergent properties (properties not present in individual motifs), depending on the initial conditions.

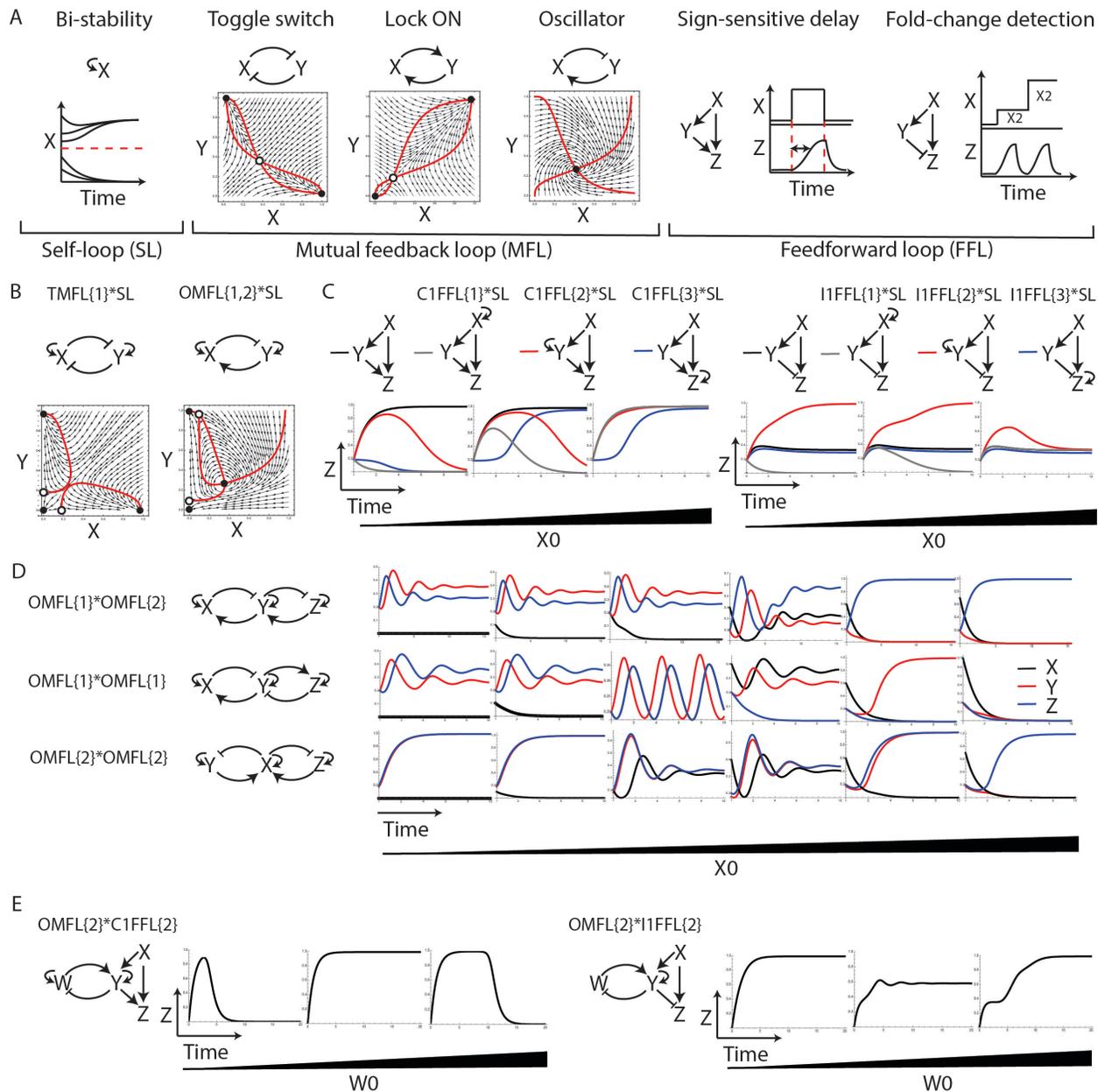

**Figure 3: Emergent dynamical properties of combinations of network motifs. A.** Three classes of circuits that we consider as minimal building blocks and examples of their dynamical properties for a given parameter values. In the phase portraits of the mutual feedback circuits the black circles represent stable fixed points, and the white circles represent unstable fixed points (models are described in the Methods section). **B.** Combinations of self-loop and mutual feedback circuits and their phase portraits (Methods). **C.** Combinations of self-loop and coherent (left panel) or incoherent (right panel) FFL and their output Z dynamical behavior. We used initial conditions of $X_0 = 0, 0.1, 0.5, Y_0 = 0.185, Z_0 = 0.19$. **D.** Three different combinations of two oscillator circuits and their dynamical properties when the initial level of X varies. We used $Y_0, Z_0 = 0.2$ throughout and $X_0 = 0, 0.1, 0.2, 0.4, 0.6, 0.7$ for each combination. **E.** Combinations of an oscillator circuit and a coherent (upper panel) or an incoherent (lower panel) FFL and their output Z dynamical behavior. We used $X_0, Z_0 = 0.01, Y_0 = 0.3$ throughout and $W_0 = 0.1, 1, 10$ for each combination.

**Emergent properties of combinations of network motifs in real networks**

We next model the dynamical behavior of several over-represented combinations of network motifs that we detected in real networks in order to exemplify the potential emergent properties that result from these combinations. The first combination of motifs that we model is the double mutual feedback motif from the *C. elegans* neuronal network (Fig. 2B). In this motif there are two pairs of neurons that mutually interact with each other (X, Y and Y, Z) while X interacts with Z only in one direction. We find that this motif is often combined in the neuronal network with other motifs of the same type where the one-direction edge from X to Z is shared between the two motifs. Considering that the interactions in the motif are all positive, each motif separately is a generalization of the Lock-ON feedback circuit where it can provide bistability such that X, Y and Z are all either turned OFF or ON. Combining two such double mutual feedback motifs links the fates of the Y and W neurons from both motifs, and can also provide a temporal order for their pulsatile behavior (Fig. 4A, Methods). We also model the same motif combination where we consider that X inhibits the activity of Z. Here, each double mutual feedback can provide oscillations for a certain range of model parameters or converge to an OFF state. When the two motifs are combined such that the X to Z inhibitory edge is shared, the oscillations of one circuit propagate to the neuron that participate in the other motif although it would not have shown oscillations autonomically with the same parameters. Moreover, it will synchronize and have the same phase as the motif it is linked with (Fig. 4B, Methods).

    The second combination of motifs that we model is of the 3-node feedback loop that was found to be enriched in the electronic circuits network (Fig. 2B). A combination of two 3-node loop circuits where one has three positive interactions (all-positive-interactions circuit) and the other has one negative and two positive interactions show emergent properties (Fig. 4C). Here the oscillations of the circuit with the negative interaction propagate to the second all-positive-interactions circuit where it oscillates with an opposite phase showing anti-phase synchronization. We note that here the emergence of oscillations in the all-positive-interactions circuit is especially interesting since on its own this circuit does not show pure oscillations even for a different choice of model parameters (Fig. 4C-D, Methods).

    These examples illustrate the complexity that can emerge from simple combinations of minimal building-block circuits in real networks.

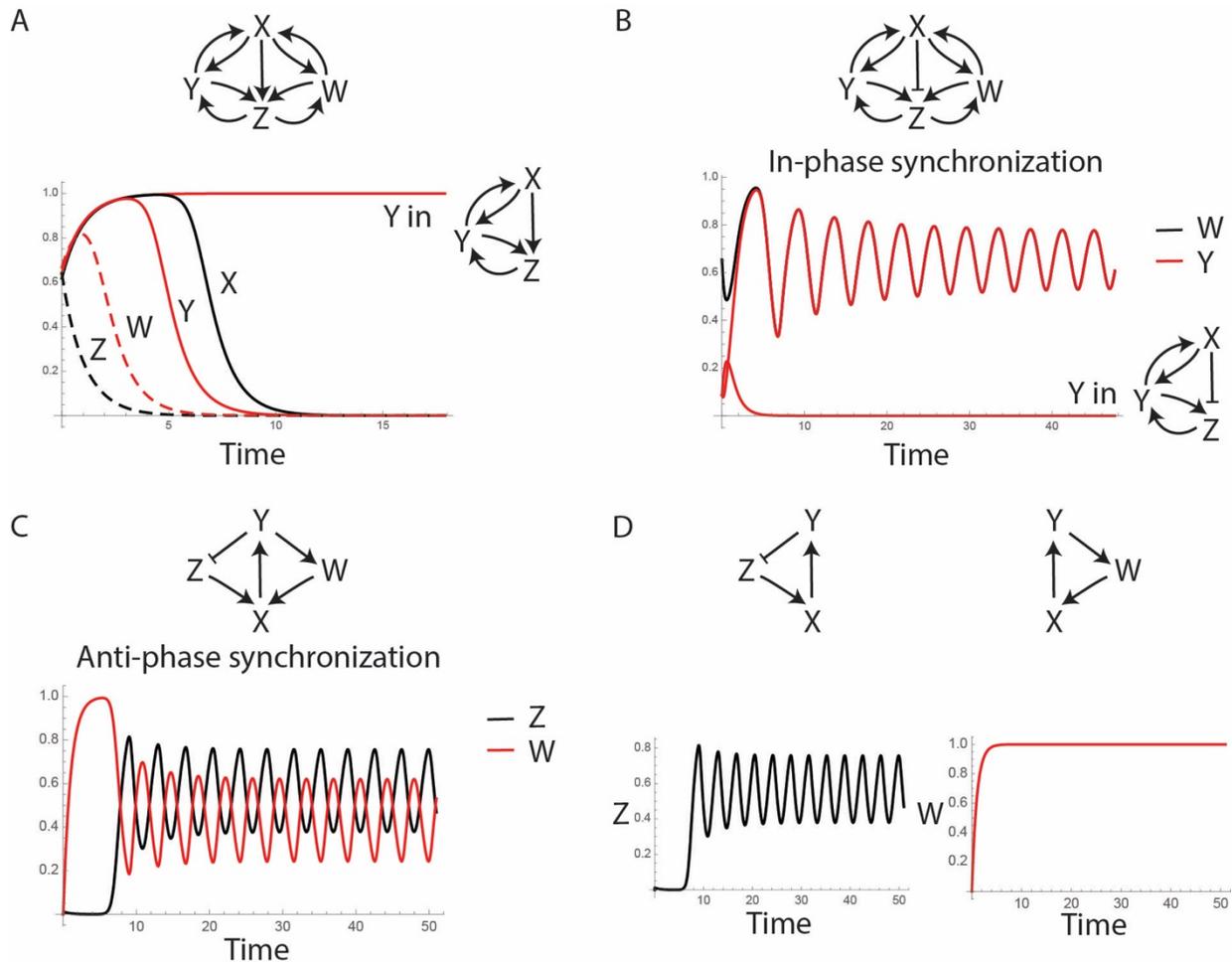

**Figure 4: Over-represented network motif combinations in real networks provide emergent properties.**
**A.** Dynamical behavior of the *C. elegans* over-represented combination of two double mutual feedback circuits where all arrows are positive interactions. Y from the separated X,Y,Z circuit with the same parameters would have converged to a high steady state level. **B.** Dynamical behavior of the *C. elegans* over-represented combination of two double mutual feedback circuits where X inhibits Z. Y from the separated X,Y,Z circuit with the same parameters would have declined to zero without oscillations. **C.** Dynamical behavior of the electronic circuits over-represented combination of two 3-node loop circuits where all edges are positive except for one where Y inhibits Z, where it shows anti-phase synchronization. **D.** W from the separated X,Y,W circuit with the same parameters (or any other choice of parameters, Methods) does not show oscillations.

**Emergent properties of interactions of network motifs**

The second way of joining two network motifs that we consider is an interaction between network motifs. Here, every motif can be considered as a hypernode in a higher-level circuit. Interaction of network motifs can thus show properties of circuits at two levels - the high-level circuit topology in which the motifs are single nodes, and the properties of the low-level circuits that are being interconnected. To illustrate this, we model an interaction of two oscillator circuits in a toggle-switch high-level topology. We find that this module shows a toggle switch between the properties of the oscillator circuits. However, the identity of the connecting nodes in each circuit and the choice of parameter values may influence the properties of the new high-level module (Fig. 5A, Methods). Similarly, two oscillator circuits that mutually activate each other show all possible combinations of their individual steady states, which is a high-level version of the Lock-ON circuit property where the motifs are either both turned OFF or ON (Fig. 5B, Methods).

Interaction between two motifs can sometimes lead to emergent properties that cannot be observed in each of the motifs separately. For example, a coherent FFL and a toggle-switch feedback circuit that mutually activate each other can exhibit oscillations which are not a property of an FFL nor a toggle-switch circuit (Methods). Joining them in this manner, creates a new path that is equivalent to an activator-inhibitor type of circuit that can generate oscillations (Fig. 5C). Oscillations can also emerge from two FFLs that interact with each other where one serves as an activator and the other as a repressor (Fig. 5D). This example demonstrates the importance of exploring the way network motifs are integrated in a large network in order to understand the potential properties that the network can exhibit. Although oscillations cannot arise from individual FFLs for any choice of parameter values (oscillations are impossible in any strictly feedforward motif), they can emerge when the FFLs mutually interact with each other in the large network.

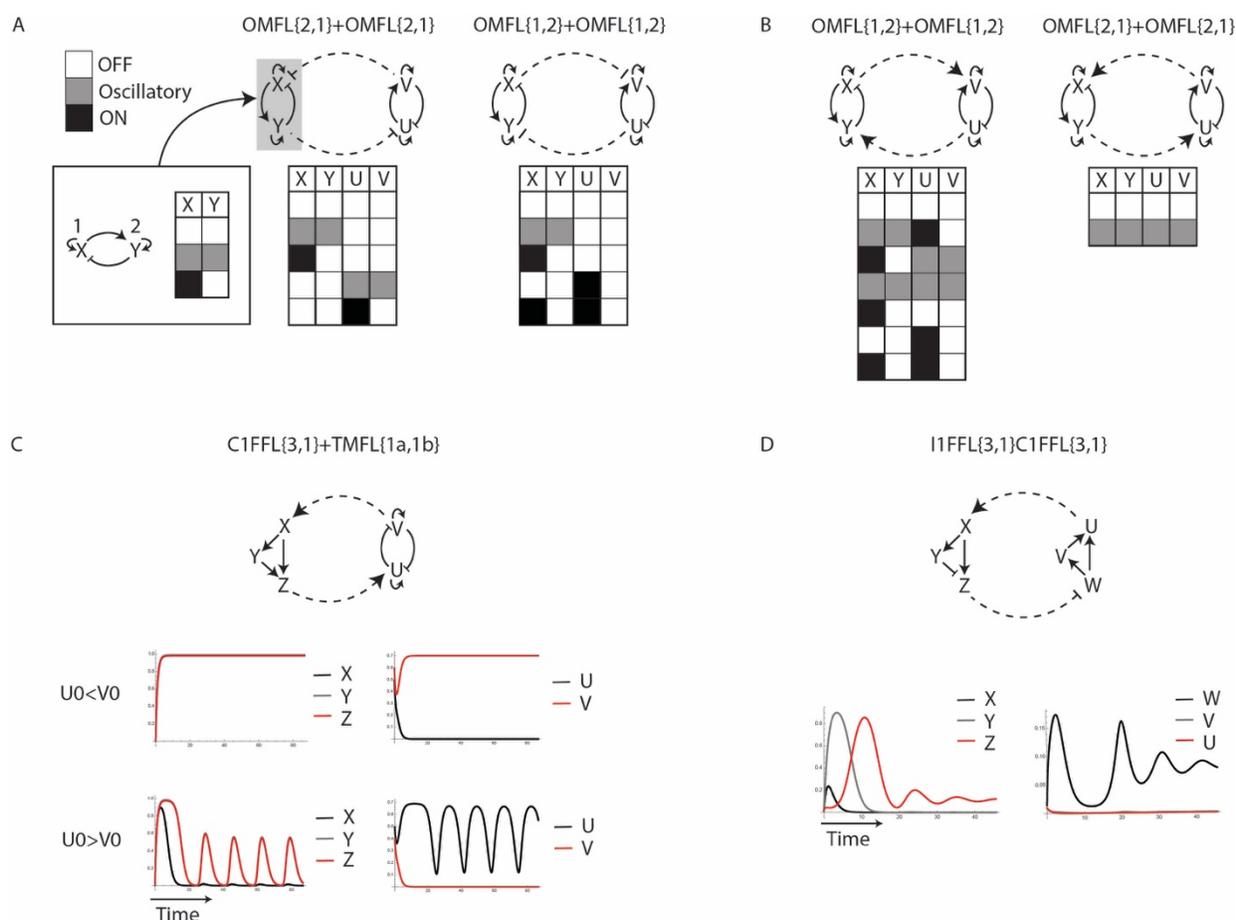

**Figure 5: Emergent dynamical properties of interactions of network motifs. A-B.** Examples of interactions of two oscillator feedback circuits where they mutually repress each other (A) and mutually activate each other (B). The stable steady states for the same choice of model parameters for all 4 circuits are shown where white squares represent the OFF state, black squares are shown when a variable is fully turned ON, and the gray squares represent damped oscillations around an intermediate level steady state. **C.** An example of a mutual activation interaction between a coherent FFL and a toggle-switch feedback circuit with the dynamic behavior of all variables for two different initial conditions with the same model parameters (see Methods). **D.** An example of an interaction between two FFLs where the coherent FFL is an activator and the incoherent FFL is a repressor, and the dynamic behavior of all variables. See Methods for the equations we used for all interactions of motifs.

**Discussion**

Complex systems are composed of multiple levels of organization. Here, we developed a theoretical framework to identify and explore the intermediate levels of organization in complex networks. We defined how two subgraphs are joined together in a network at two levels by sharing at least one vertex or by being directly linked by at least one edge. We developed a new method to reveal how network motifs are assembled into hyper-motifs in real complex networks and demonstrated it in diverse evolved and designed networks. Finally, we used a nonlinear class of models to explore emergent properties of combinations and interactions of network motifs including the autoregulation, feedforward loop and several feedback circuits.

Exploring properties of combinations of building-block circuits provides a way to rigorously explore emergence in complex systems and to define new levels of organization based on functional modules that provide important and emergent properties. Applying our framework on real networks allowed us to reveal patterns in the inner topological structure of the network and to make specific predictions on the behavior that is expected to emerge at the mesoscale level of the network.

The framework presented here could be useful for developing a new mnemonic way to categorize the nodes in a network into distinct classes based on their position in the network's motifs and their higher-order modules. This categorization is complement to previously suggested criteria including the nodes degree and contrabillity properties (*43*) where nodes with the same level of connectivity can play different roles in the network's modules. Identifying these categories of nodes in real networks can highlight the key drivers of the network's emergent properties. For example, in biology this could have important implications in pinpointing genes that are key for essential biological processes as well as genes with a potential to drive a biological system into a pathological state.

**Acknowledgments:**

We thank Drs. Uri Alon, Yang-Yu Liu, Shlomo Havlin, Andre Levchenko, Baruch Barzel, Yuval Hart, Avi Mayo, Mor Nitzan and members of the Medzhitov lab for discussions. MA was supported by the EMBO Long-Term Fellowship (ALTF 304-2019) and the Zuckerman STEM Leadership program. RM is supported by the Howard Hughes Medical Institute and the Blavatnik Family Foundation.

**Author contributions:**

Conceptualization: MA, RM

Methodology and Formal Analysis: MA

Writing: MA, RM

**Competing interests:** Authors declare that they have no competing interests.

**Materials and Methods**

**Counting all possible combinations and interactions of motifs in our framework**

Considering two network motifs, A and B, with $n_A$ and $n_B$ nodes, respectively. A and B can be *combined* by sharing $N_v$ nodes such that the following condition applies: $1 \leq N_v < min(n_A, n_B)$. This condition ensures that the topology of each motif is kept as a subgraph of the combined circuit. The identity of the $N_v$ nodes that are shared among the two motifs define a core topology of a combination of the two motifs. Each core topology can be extended to include additional edges between pairs of nodes that do not participate in the same motif (fig. S1A).

A and B *interact* when they have at least one edge that directly links them. Since each node from motif A can interact with each node from motif B, the maximal number of possible edges is the number of pairs of nodes that don't participate in the same motif, $n_A n_B$. For directed networks where every pair of nodes can have two directed edges this number should be multiplied by 2. The number of possible topologies of interaction between motifs A and B is therefore $2^{2n_A n_B}$ for directed networks and $2^{n_A n_B}$ for undirected networks (fig. S1B). We note that in both combinations and interactions of motifs, there could be circuits that are isomorphic to each other and therefore the number of unique topologies may be smaller than the maximal number of circuit topologies.

**Detecting enriched combinations of network motifs in real networks**

Consider a network G with N nodes and E edges. To detect over- and under-represented combinations of network motifs we followed the following steps:

1. Identify network motifs of up to 3-nodes (see remark below). We used MFinder (*7*) and also verified consistency with finding network motifs in Mathematica using the IGraphM package and its function IGRewire to randomize the network, and IGMotifs to find the motifs and their frequencies. We denote $N_m$ as the set of nodes in G that participate in the network motifs.
2. Categorize nodes in $N_m$ to different groups according to their roles in the network motifs: $\{n_1, n_2, \ldots, n_k\}$, where k is the number of the different network motifs roles in G. Note that certain nodes in $N_m$ may appear in more than one of the $\{n_1, n_2, \ldots, n_k\}$ groups (which are the cases where network motifs are combined).
3. Compute the Jaccard index for all $k(k-1)/2$ pairs of motif roles: $J(n_i, n_j) = |n_i \cap n_j|/|n_i \cup n_j|$ which is the size of the intersection of nodes in roles i,j divided by the size of their union.
4. For nodes that appear more than once in the same role, we compute the Jaccard index by the ratio of the number of nodes that appear in a network motif role more than once to the total number of nodes that participate in that motif role.
5. We use the MFinder package to create 100 random networks that have N nodes, E edges, the same incoming and outgoing edges per node and the same frequency of all subgraphs up to 3-nodes subgraphs.
6. We repeat steps 2-4 for all the random networks, and compute $J_{rand}(n_i, n_j)$ for all pairs of motif roles for the random networks.
7. For every i,j such that $i, j \in \{1, \ldots, k\}$ we compute the Z-score of $J(n_i, n_j)$ of network G from the distribution of $\{J_{rand}(n_i, n_j)\}$: $Z_{ij} = (J(n_i, n_j) - mean(\{J_{rand}(n_i, n_j)\}))/std(\{J_{rand}(n_i, n_j)\}$ and compute the p-value by estimating the cumulative density (CDF) of $Z_{ij}$ for a normal distribution with zero mean and unit variance. We use the function NormalPValue in the package HypothesisTesting in Mathematica to compute the p-value. We then use the Binjamini-Hochberg procedure to correct for multiple hypothesis testing which provides us a corrected q-value for each pair of motif roles.
8. We consider over- and under-represented motif combinations if their q-value is smaller than 0.05, where $Z_{ij} > 0$ for over-represented combinations and $Z_{ij} < 0$ for under-represented combinations.

Our method detects situations in which two network motifs are joined by sharing certain nodes. The difference between our method and detecting network motifs of a larger size is that our method does not detect a specific subgraph as a recurring pattern, but rather finds classes of subgraphs in which the two network motifs in question are combined in a certain way. This means that in certain networks the exact topology of the combination of the two network motifs may include additional edges that are in accordance with our definition of a core topology of a combination of network motifs and its possible extensions (see fig. S1A). In figure S2A we show the frequencies of the core topologies of each over-represented motif combination and its extensions. We find that in most cases the most enriched motif combinations are with a core topology (without additional edges besides the edges of the network motifs).

To analyze the word adjacency network, we downsampled the network (G) using the following steps:
1. Define the size of the downsampling, $sz$, as a parameter that can be tuned.
2. We will construct a list of nodes which we sample, $s$, from the large network G.
3. Randomly choose one node from G, $s_0$, to be the first entry of $s$.
4. Randomly choose a node that is a member of the neighborhood of $s_0$, $s_1 \in N_{s_0}$, to be the next entry of $s$.
5. Insert additional entries to $s$ for $i \in \{2, \ldots, sz/3\}$ by randomly choosing a node that is a member of the neighborhood of $s_{i-1}$ with probability of 85%, or a node that is a member of the neighborhood of $s_0$.
6. Certain entries of $s$ could be repeated more than once. If the length of the unique list of $s$ is smaller than $(sz/3)/2$ - consider a new $s_0$ which is a randomly chosen node from G. Otherwise, keep the previous $s_0$.
7. Insert additional entries to $s$ for $i \in \{sz/3 + 1, \ldots, 2sz/3\}$ by randomly choosing a node that is a member of the neighborhood of $s_{i-1}$ with probability of 85%, or a node that is a member of the neighborhood of $s_0$.
8. If the length of the unique list of $s$ is smaller than $(sz/3)$ - consider a new $s_0$ which is a randomly chosen node from G. Otherwise, keep the previous $s_0$.
9. Insert additional entries to $s$ for $i \in \{2sz/3 + 1, \ldots, sz\}$ by randomly choosing a node that is a member of the neighborhood of $s_{i-1}$ with probability of 85%, or a node that is a member of the neighborhood of $s_0$.
10. The downsampled network $G_d$ is the network that contains nodes from the unique list of $s$.

We repeated the downsampling procedure multiple times with varying sampling size ($sz$), where we check that the downsampled networks show a similar degree distribution and the same network motifs as the complete large network.

**Modeling known network motifs and their combinations**

We model several previously identified network motifs where we consider that each node has a linear removal term and that interactions between nodes can be described by Hill functions. A positive interaction from node X to node Y is modeled by an increasing Hill function, $Y^{n_{xy}}/(k_{xy}{}^{n_{xy}} + Y^{n_{xy}})$, and a negative interaction is modeled by a decreasing Hill function, $k_{xy}{}^{n_{xy}}/(k_{xy}{}^{n_{xy}} + Y^{n_{xy}})$. Each interaction is therefore described by two parameters: $n_{xy}$ is the cooperativity coefficient, and $k_{xy}$ is the level at which the interaction effect reaches halfway of its maximal level.

We now list the models and parameter values that we used in Figure 3:
- Positive self-loop:
  M3. $X' = X^{n_{xx}}/(k_{xx}{}^{n_{xx}} + X^{n_{xx}}) - X$

The cooperativity, $n_{xx}$, has to be at least 2 for this model of a self-loop circuit to show bistability.
- Toggle-switch feedback circuit:
  M4. $X' = k_{yx}{}^{n_{yx}}/(k_{yx}{}^{n_{yx}} + Y^{n_{yx}}) - X$
  M5. $Y' = k_{xy}{}^{n_{xy}}/(k_{xy}{}^{n_{xy}} + X^{n_{xy}}) - Y$
- Lock-ON feedback circuit:
  M6. $X' = Y^{n_{yx}}/(k_{yx}{}^{n_{yx}} + Y^{n_{yx}}) - X$
  M7. $Y' = X^{n_{xy}}/(k_{xy}{}^{n_{xy}} + X^{n_{xy}}) - Y$
- Oscillator feedback circuit:
  M8. $X' = Y^{n_{yx}}/(k_{yx}{}^{n_{yx}} + Y^{n_{yx}}) - X$
  M9. $Y' = k_{xy}{}^{n_{xy}}/(k_{xy}{}^{n_{xy}} + X^{n_{xy}}) - Y$

The phase portraits for the feedback circuits in Figure 3A are plotted using the above models (eqs. M4-M9) with $n_{yx} = n_{xy} = 3, k_{yx} = k_{xy} = 0.3$. We used the functions Streamplot and Contourplot in Mathematica 12.1.1.0 to plot the phase portraits of the circuits.
- Combination of positive self-loops and a toggle-switch circuit:
  M10.     $X' = (X^{n_{xx}}/(k_{xx}{}^{n_{xx}} + X^{n_{xx}}))(k_{yx}{}^{n_{yx}}/(k_{yx}{}^{n_{yx}} + Y^{n_{yx}})) - X$, $n_{xx} = n_{yx} = 3, k_{xx} = k_{yx} = 0.3$
  M11.     $Y' = (Y^{n_{yy}}/(k_{yy}{}^{n_{yy}} + Y^{n_{yy}}))(k_{xy}{}^{n_{xy}}/(k_{xy}{}^{n_{xy}} + X^{n_{xy}})) - Y$, $n_{yy} = n_{xy} = 3, k_{yy} = k_{xy} = 0.3$
- Combination of positive self-loops and an oscillator circuit:
  M12. $X' = (X^{n_{xx}}/(k_{xx}{}^{n_{xx}} + X^{n_{xx}}))(Y^{n_{yx}}/(k_{yx}{}^{n_{yx}} + Y^{n_{yx}})) - X$, $n_{xx} = n_{yx} = 3, k_{xx} = 0.2, k_{yx} = 0.3$
  M13.     $Y' = (Y^{n_{yy}}/(k_{yy}{}^{n_{yy}} + Y^{n_{yy}}))(k_{xy}{}^{n_{xy}}/(k_{xy}{}^{n_{xy}} + X^{n_{xy}})) - Y$, $n_{yy} = n_{xy} = 3, k_{yy} = 0.2, k_{xy} = 0.3$
- Combination of a positive self-loop and feedforward loop (FFL) circuits:
  The model for a coherent FFL:
  M14. $X' = 1 - X$
  M15. $Y' = X^{n_{xy}}/(k_{xy}{}^{n_{xy}} + X^{n_{xy}}) - Y$, $n_{xy} = 1, k_{xy} = 0.01$
  M16. $Z' = (X^{n_{xz}}/(k_{xz}{}^{n_{xz}} + X^{n_{xz}}))(Y^{n_{yz}}/(k_{yz}{}^{n_{yz}} + Y^{n_{yz}})) - Z$, $n_{xz} = n_{yz} = 1, k_{xz} = k_{yz} = 0.01$

  The model for an incoherent FFL:
  M17. $X' = 1 - X$
  M18. $Y' = X^{n_{xy}}/(k_{xy}{}^{n_{xy}} + X^{n_{xy}}) - Y$, $n_{xy} = 1, k_{xy} = 0.01$
  M19. $Z' = (X^{n_{xz}}/(k_{xz}{}^{n_{xz}} + X^{n_{xz}}))(k_{yz}{}^{n_{yz}}/(k_{yz}{}^{n_{yz}} + Y^{n_{yz}})) - Z$, $n_{xz} = n_{yz} = 1, k_{xz} = 0.01, k_{yz} = 0.5$

To combine the FFL circuits with a self-loop on one of its nodes, i, we multiply the production term of variable i by the following term: $(i^{n_{ii}}/(k_{ii}{}^{n_{ii}} + i^{n_{ii}}))$, $n_{ii} = 3, k_{ii} = 0.3$ in the coherent

FFL and $k_{xx} = k_{yy} = 0.3, k_{zz} = 0.15$ in the incoherent FFL. In the SI we explore combinations of FFL circuits with self loops on the different FFL nodes where we assume that X is not a dynamical variable but rather rises in a step function manner. This assumption allows us to explore the phase portraits of Y and Z and to compare their nullclines when the self loop appears on Y or Z (fig. S3A-B).

- Combinations of two oscillator circuits:
  First combination:
  M20. $X' = (X^{n_{xx}}/(k_{xx}^{n_{xx}} + X^{n_{xx}}))(Y^{n_{yx}}/(k_{yx}^{n_{yx}} + Y^{n_{yx}})) - X$, $n_{xx} = n_{yx} = 3, k_{xx} = 0.2, k_{yx} = 0.3$

  M21. $Y' = (Y^{n_{yy}}/(k_{yy}^{n_{yy}} + Y^{n_{yy}}))(k_{xy}^{n_{xy}}/(k_{xy}^{n_{xy}} + X^{n_{xy}}))(Z^{n_{zy}}/(k_{zy}^{n_{zy}} + Z^{n_{zy}})) - Y$, $n_{yy} = n_{xy} = n_{zy} = 3, k_{yy} = 0.2, k_{xy} = k_{zy} = 0.3$

  M22. $Z' = (Z^{n_{zz}}/(k_{zz}^{n_{zz}} + Z^{n_{zz}}))(k_{yz}^{n_{yz}}/(k_{yz}^{n_{yz}} + Y^{n_{yz}})) - Z$, $n_{zz} = n_{yz} = 3, k_{zz} = 0.2, k_{yz} = 0.3$

  Second combination:
  M22. $X' = (X^{n_{xx}}/(k_{xx}^{n_{xx}} + X^{n_{xx}}))(Y^{n_{yx}}/(k_{yx}^{n_{yx}} + Y^{n_{yx}})) - X$, $n_{xx} = n_{yx} = 3, k_{xx} = 0.2, k_{yx} = 0.3$

  M23. $Y' = (Y^{n_{yy}}/(k_{yy}^{n_{yy}} + Y^{n_{yy}}))(k_{xy}^{n_{xy}}/(k_{xy}^{n_{xy}} + X^{n_{xy}}))(k_{zy}^{n_{zy}}/(k_{zy}^{n_{zy}} + Z^{n_{zy}})) - Y$
  $n_{yy} = n_{xy} = n_{zy} = 3, k_{yy} = 0.2, k_{xy} = k_{zy} = 0.3$

  M24. $Z' = (Z^{n_{zz}}/(k_{zz}^{n_{zz}} + Z^{n_{zz}}))(Y^{n_{yz}}/(k_{yz}^{n_{yz}} + Y^{n_{yz}})) - Z$, $n_{zz} = n_{yz} = 3, k_{zz} = 0.2, k_{yz} = 0.3$

  Third combination:
  M25. $X' = (X^{n_{xx}}/(k_{xx}^{n_{xx}} + X^{n_{xx}}))(k_{yx}^{n_{yx}}/(k_{yx}^{n_{yx}} + Y^{n_{yx}})) - X$, $n_{xx} = n_{yx} = 3, k_{xx} = 0.2, k_{yx} = 0.3$

  M26. $Y' = (Y^{n_{yy}}/(k_{yy}^{n_{yy}} + Y^{n_{yy}}))(X^{n_{xy}}/(k_{xy}^{n_{xy}} + X^{n_{xy}}))(Z^{n_{zy}}/(k_{zy}^{n_{zy}} + Z^{n_{zy}})) - Y$

  $n_{yy} = n_{xy} = n_{zy} = 3, k_{yy} = 0.2, k_{xy} = k_{zy} = 0.3$

  M27. $Z' = (Z^{n_{zz}}/(k_{zz}^{n_{zz}} + Z^{n_{zz}}))(k_{yz}^{n_{yz}}/(k_{yz}^{n_{yz}} + Y^{n_{yz}})) - Z$, $n_{zz} = n_{yz} = 3, k_{zz} = 0.2, k_{yz} = 0.3$

- Combination of an oscillator circuit and a coherent FFL through the FFL's intermediate node:
  M27. $X' = 1 - X$

  M28. $Y' = (Y^{n_{yy}}/(k_{yy}^{n_{yy}} + Y^{n_{yy}}))(X^{n_{xy}}/(k_{xy}^{n_{xy}} + X^{n_{xy}}))(W^{n_{wy}}/(k_{wy}^{n_{wy}} + W^{n_{wy}})) - Y$

  $n_{yy} = n_{wy} = 3, n_{xy} = 1, k_{yy} = 0.2, k_{wy} = 0.3, k_{xy} = 0.01$

  M29. $Z' = (X^{n_{xz}}/(k_{xz}^{n_{xz}} + X^{n_{xz}}))(Y^{n_{yz}}/(k_{yz}^{n_{yz}} + Y^{n_{yz}})) - Z$, $n_{xz} = n_{yz} = 1, k_{xz} = k_{yz} = 0.01$

  M30. $W' = (W^{n_{ww}}/(k_{ww}^{n_{ww}} + W^{n_{ww}}))(k_{yw}^{n_{yw}}/(k_{yw}^{n_{yw}} + Y^{n_{yw}})) - W$, $n_{ww} = n_{yw} = 3, k_{ww} = 0.2, k_{yw} = 0.3$

- Combination of an oscillator circuit and an incoherent FFL through the FFL's intermediate node:
  M27. $X' = 1 - X$

  M28. $Y' = (Y^{n_{yy}}/(k_{yy}^{n_{yy}} + Y^{n_{yy}}))(X^{n_{xy}}/(k_{xy}^{n_{xy}} + X^{n_{xy}}))(W^{n_{wy}}/(k_{wy}^{n_{wy}} + W^{n_{wy}})) - Y$

$$n_{yy} = n_{wy} = 3, n_{xy} = 1, k_{yy} = 0.2, k_{wy} = 0.3, k_{xy} = 0.01$$

M29. $Z' = (X^{n_{xz}}/(k_{xz}^{n_{xz}} + X^{n_{xz}}))(k_{yz}^{n_{yz}}/(k_{yz}^{n_{yz}} + Y^{n_{yz}})) - Z$, $n_{xz} = n_{yz} = 1, k_{xz} = 0.01, k_{yz} = 0.5$

M30. $W' = (W^{n_{ww}}/(k_{ww}^{n_{ww}} + W^{n_{ww}}))(k_{yw}^{n_{yw}}/(k_{yw}^{n_{yw}} + Y^{n_{yw}})) - W$, $n_{ww} = n_{yw} = 3, k_{ww} = 0.2, k_{yw} = 0.3$

In the SI we plot the dynamical behavior of all possible combinations between these network motifs for a given choice of model parameters (fig. S3C-R).

**Modeling known network motifs and their interactions**

We now list the models and parameter values that we used to model interactions of network motifs that are shown in Figure 5:

- Interaction between two oscillator circuits in a toggle switch topology (option #1):

    M31. $X' = (X^{n_{xx}}/(k_{xx}^{n_{xx}} + X^{n_{xx}}))(k_{yx}^{n_{yx}}/(k_{yx}^{n_{yx}} + Y^{n_{yx}}))(k_{vx}^{n_{vx}}/(k_{vx}^{n_{vx}} + V^{n_{vx}})) - X$

    M32. $Y' = (Y^{n_{yy}}/(k_{yy}^{n_{yy}} + Y^{n_{yy}}))(X^{n_{xy}}/(k_{xy}^{n_{xy}} + X^{n_{xy}})) - Y$

    M33. $U' = (U^{n_{uu}}/(k_{uu}^{n_{uu}} + U^{n_{uu}}))(k_{yu}^{n_{yu}}/(k_{yu}^{n_{yu}} + Y^{n_{yu}}))(k_{vu}^{n_{vu}}/(k_{vu}^{n_{vu}} + V^{n_{vu}})) - U$

    M34. $V' = (V^{n_{vv}}/(k_{vv}^{n_{vv}} + V^{n_{vv}}))(U^{n_{uv}}/(k_{uv}^{n_{uv}} + U^{n_{uv}})) - V$

    Where we used the following parameter values:
    $k_{ii} = 0.2, n_{ii} = 3, k_{xy} = k_{yx} = k_{uv} = k_{vu} = 0.3, k_{yu} = k_{vx} = 0.01, n_{xy} = n_{yx} = n_{uv} = n_{vu} = 3, n_{yu} = n_{vx} = 1$

- Interaction between two oscillator circuits in a toggle switch topology (option #2):

    M35. $X' = (X^{n_{xx}}/(k_{xx}^{n_{xx}} + X^{n_{xx}}))(k_{yx}^{n_{yx}}/(k_{yx}^{n_{yx}} + Y^{n_{yx}}))) - X$

    M36. $Y' = (Y^{n_{yy}}/(k_{yy}^{n_{yy}} + Y^{n_{yy}}))(X^{n_{xy}}/(k_{xy}^{n_{xy}} + X^{n_{xy}}))(k_{uy}^{n_{uy}}/(k_{uy}^{n_{uy}} + U^{n_{uy}})) - Y$

    M37. $U' = (U^{n_{uu}}/(k_{uu}^{n_{uu}} + U^{n_{uu}}))(k_{vu}^{n_{vu}}/(k_{vu}^{n_{vu}} + V^{n_{vu}})) - U$

    M38. $V' = (V^{n_{vv}}/(k_{vv}^{n_{vv}} + V^{n_{vv}}))(U^{n_{uv}}/(k_{uv}^{n_{uv}} + U^{n_{uv}}))(k_{xv}^{n_{xv}}/(k_{xv}^{n_{xv}} + X^{n_{xv}}) - V$

    Where we used the following parameter values:
    $k_{ii} = 0.2, n_{ii} = 3, k_{xy} = k_{yx} = k_{uv} = k_{vu} = 0.3, k_{uy} = k_{xv} = 0.01, n_{xy} = n_{yx} = n_{uv} = n_{vu} = 3, n_{uy} = n_{xv} = 1$

- Interaction between two oscillator circuits in a lock-ON topology (option #1):

    M39. $X' = (X^{n_{xx}}/(k_{xx}^{n_{xx}} + X^{n_{xx}}))(k_{yx}^{n_{yx}}/(k_{yx}^{n_{yx}} + Y^{n_{yx}}))) - X$

    M40. $Y' = (Y^{n_{yy}}/(k_{yy}^{n_{yy}} + Y^{n_{yy}}))(X^{n_{xy}}/(k_{xy}^{n_{xy}} + X^{n_{xy}}))(U^{n_{uy}}/(k_{uy}^{n_{uy}} + U^{n_{uy}})) - Y$

    M41. $U' = (U^{n_{uu}}/(k_{uu}^{n_{uu}} + U^{n_{uu}}))(k_{vu}^{n_{vu}}/(k_{vu}^{n_{vu}} + V^{n_{vu}})) - U$

    M42. $V' = (V^{n_{vv}}/(k_{vv}^{n_{vv}} + V^{n_{vv}}))(U^{n_{uv}}/(k_{uv}^{n_{uv}} + U^{n_{uv}}))(X^{n_{xv}}/(k_{xv}^{n_{xv}} + X^{n_{xv}}) - V$

    Where we used the following parameter values:
    $k_{ii} = 0.2, n_{ii} = 3, k_{xy} = k_{yx} = k_{uv} = k_{vu} = 0.3, k_{yu} = k_{vx} = 0.01, n_{xy} = n_{yx} = n_{uv} = n_{vu} = 3, n_{yu} = n_{vx} = 1$

- Interaction between two oscillator circuits in a lock-ON topology (option #2):

    M43. $X' = (X^{n_{xx}}/(k_{xx}^{n_{xx}} + X^{n_{xx}}))(k_{yx}^{n_{yx}}/(k_{yx}^{n_{yx}} + Y^{n_{yx}}))(V^{n_{vx}}/(k_{vx}^{n_{vx}} + V^{n_{vx}})) - X$

    M44. $Y' = (Y^{n_{yy}}/(k_{yy}^{n_{yy}} + Y^{n_{yy}}))(X^{n_{xy}}/(k_{xy}^{n_{xy}} + X^{n_{xy}})) - Y$

M45. $U' = (U^{n_{uu}}/(k_{uu}{}^{n_{uu}} + U^{n_{uu}}))(Y^{n_{yu}}/(k_{yu}{}^{n_{yu}} + Y^{n_{yu}}))(k_{vu}{}^{n_{vu}}/(k_{vu}{}^{n_{vu}} + V^{n_{vu}})) - U$

M46. $V' = (V^{n_{vv}}/(k_{vv}{}^{n_{vv}} + V^{n_{vv}}))(U^{n_{uv}}/(k_{uv}{}^{n_{uv}} + U^{n_{uv}})) - V$

Where we used the following parameter values:
$k_{ii} = 0.2, n_{ii} = 3, k_{xy} = k_{yx} = k_{uv} = k_{vu} = 0.3, k_{uy} = k_{xv} = 0.01, n_{xy} = n_{yx} = n_{uv} = n_{vu} = 3, n_{uy} = n_{xv} = 1$

- Interaction between a coherent FFL and a toggle switch circuit:

    M47. $X' = (V^{n_{vx}}/(k_{vx}{}^{n_{vx}} + V^{n_{vx}})) - X$

    M48. $Y' = (X^{n_{xy}}/(k_{xy}{}^{n_{xy}} + X^{n_{xy}})) - Y$

    M49. $Z' = (Y^{n_{yz}}/(k_{yz}{}^{n_{yz}} + Y^{n_{yz}}))(X^{n_{xz}}/(k_{xz}{}^{n_{xz}} + X^{n_{xz}})) - Z$

    M50. $U' = (U^{n_{uu}}/(k_{uu}{}^{n_{uu}} + U^{n_{uu}}))(Z^{n_{zu}}/(k_{zu}{}^{n_{zu}} + Z^{n_{zu}}))(k_{vu}{}^{n_{vu}}/(k_{vu}{}^{n_{vu}} + V^{n_{vu}})) - U$

    M51. $V' = (V^{n_{vv}}/(k_{vv}{}^{n_{vv}} + V^{n_{vv}}))(k_{uv}{}^{n_{uv}}/(k_{uv}{}^{n_{uv}} + U^{n_{uv}})) - V$

    Where we used the following parameter values:
    $k_{ii} = 0.3, n_{ii} = 1, k_{uv} = k_{vu} = 0.3, k_{vx} = k_{xy} = k_{xz} = k_{yz} = k_{zu} = 0.01, n_{uv} = n_{vu} = 3, n_{vx} = n_{xy} = n_{xz} = n_{yz} = n_{zu} = 1$

- Interaction between a coherent and an incoherent FFLs:

    M52. $X' = (U^{n_{ux}}/(k_{ux}{}^{n_{ux}} + U^{n_{ux}})) - X$

    M53. $Y' = (X^{n_{xy}}/(k_{xy}{}^{n_{xy}} + X^{n_{xy}})) - Y$

    M54. $Z' = (k_{yz}{}^{n_{yz}}/(k_{yz}{}^{n_{yz}} + Y^{n_{yz}}))(X^{n_{xz}}/(k_{xz}{}^{n_{xz}} + X^{n_{xz}})) - Z$

    M55. $W' = (k_{zw}{}^{n_{zw}}/(k_{zw}{}^{n_{zw}} + Z^{n_{zw}})) - W$

    M56. $V' = (W^{n_{wv}}/(k_{wv}{}^{n_{wv}} + W^{n_{wv}})) - V$

    M57. $U' = (V^{n_{vu}}/(k_{vu}{}^{n_{vu}} + V^{n_{vu}}))(W^{n_{wu}}/(k_{wu}{}^{n_{wu}} + W^{n_{wu}})) - U$

    Where we used the following parameter values:
    $k_{ux} = k_{xy} = k_{xz} = k_{yz} = k_{zw} = 0.01, k_{vu} = 13.4, k_{wu} = 55.9, k_{wv} = 70.8, n_{ij} = 1$

In the last example of an interaction between two FFLs, the oscillations are an emergent property since they are absent from FFL circuits on their own. To prove that, consider the Jacobian matrix of the FFL circuit: $J = [-1\ 0\ 0]\backslash[1\ -1\ 0]\backslash[1\ 1\ -1]$ which is a lower triangular matrix. The eigenvalues of a lower (or an upper) triangular matrix are the entries on its diagonal and therefore they are always -1 in our case and cannot be complex.

**Modeling enriched combinations in real networks**

We now list the models and parameter values that we used to model several observed over-represented combinations of network motifs that are shown in Figure 2:

- Combination of two double mutual feedback circuits where all interactions are positive:

    M58. $X' = (W^{n_{wx}}/(k_{wx}{}^{n_{wx}} + W^{n_{wx}}))(Y^{n_{yx}}/(k_{yx}{}^{n_{yx}} + Y^{n_{yx}})) - X$

    M59. $Y' = (Z^{n_{zy}}/(k_{zy}{}^{n_{zy}} + Z^{n_{zy}}))(X^{n_{xy}}/(k_{xy}{}^{n_{xy}} + X^{n_{xy}})) - Y$

    M60. $Z' = (X^{n_{xz}}/(k_{xz}{}^{n_{xz}} + X^{n_{xz}}))(Y^{n_{yz}}/(k_{yz}{}^{n_{yz}} + Y^{n_{yz}}))(W^{n_{wz}}/(k_{wz}{}^{n_{wz}} + W^{n_{wz}})) - Z$

    M61. $W' = (X^{n_{xw}}/(k_{xw}{}^{n_{xw}} + X^{n_{xw}}))(Z^{n_{zw}}/(k_{zw}{}^{n_{zw}} + Z^{n_{zw}})) - W$

    Where we used the following parameter values:
    $k_{wx} = k_{xw} = k_{xy} = k_{yx} = k_{xz} = k_{yz} = k_{zy} = 0.01, k_{wz} = 3.38, k_{zw} = 0.16, n_{ij} = 3$

- Combination of two double mutual feedback circuits where all interactions are positive except for the edge between X and Z where X inhibits Z:

M62. $X' = (W^{n_{wx}}/(k_{wx}^{n_{wx}} + W^{n_{wx}}))(Y^{n_{yx}}/(k_{yx}^{n_{yx}} + Y^{n_{yx}})) - X$

M63. $Y' = (Z^{n_{zy}}/(k_{zy}^{n_{zy}} + Z^{n_{zy}}))(X^{n_{xy}}/(k_{xy}^{n_{xy}} + X^{n_{xy}})) - Y$

M64. $Z' = (k_{xz}^{n_{xz}}/(k_{xz}^{n_{xz}} + X^{n_{xz}}))(Y^{n_{yz}}/(k_{yz}^{n_{yz}} + Y^{n_{yz}}))(W^{n_{wz}}/(k_{wz}^{n_{wz}} + W^{n_{wz}})) - Z$

M65. $W' = (X^{n_{xw}}/(k_{xw}^{n_{xw}} + X^{n_{xw}}))(Z^{n_{zw}}/(k_{zw}^{n_{zw}} + Z^{n_{zw}})) - W$

Where we used the following parameter values:

$k_{wx} = k_{xw} = k_{xy} = k_{xz} = k_{yz} = k_{zy} = k_{wz} = k_{zw} = 0.01, k_{yx} = 1.87, n_{ij} = 3$

- Combination of two 3-node loop circuits where all interactions are positive except for the edge between Y and Z where Y inhibits Z:

M66. $X' = (W^{n_{wx}}/(k_{wx}^{n_{wx}} + W^{n_{wx}}))(Z^{n_{zx}}/(k_{zx}^{n_{zx}} + Z^{n_{zx}})) - X$

M67. $Y' = (X^{n_{xy}}/(k_{xy}^{n_{xy}} + X^{n_{xy}})) - Y$

M68. $Z' = (k_{yz}^{n_{yz}}/(k_{yz}^{n_{yz}} + Y^{n_{yz}})) - Z$

M69. $W' = (Y^{n_{yw}}/(k_{yw}^{n_{yw}} + Y^{n_{yw}})) - W$

Where we used the following parameter values:

$k_{wx} = k_{xy} = k_{yw} = k_{yz} = 0.01, k_{zx} = 4.7, n_{ij} = 3$

In the last example, undamped oscillations emerge in the all-positive-interactions 3-node loop circuit. The reason is that pure oscillations require a negative feedback and a delay. However, we note that this all-positive-interactions circuit can show damped oscillations (with a spiral fixed point) which in the presence of noise can become undamped oscillations (Alon 2019).

**SI Figures**

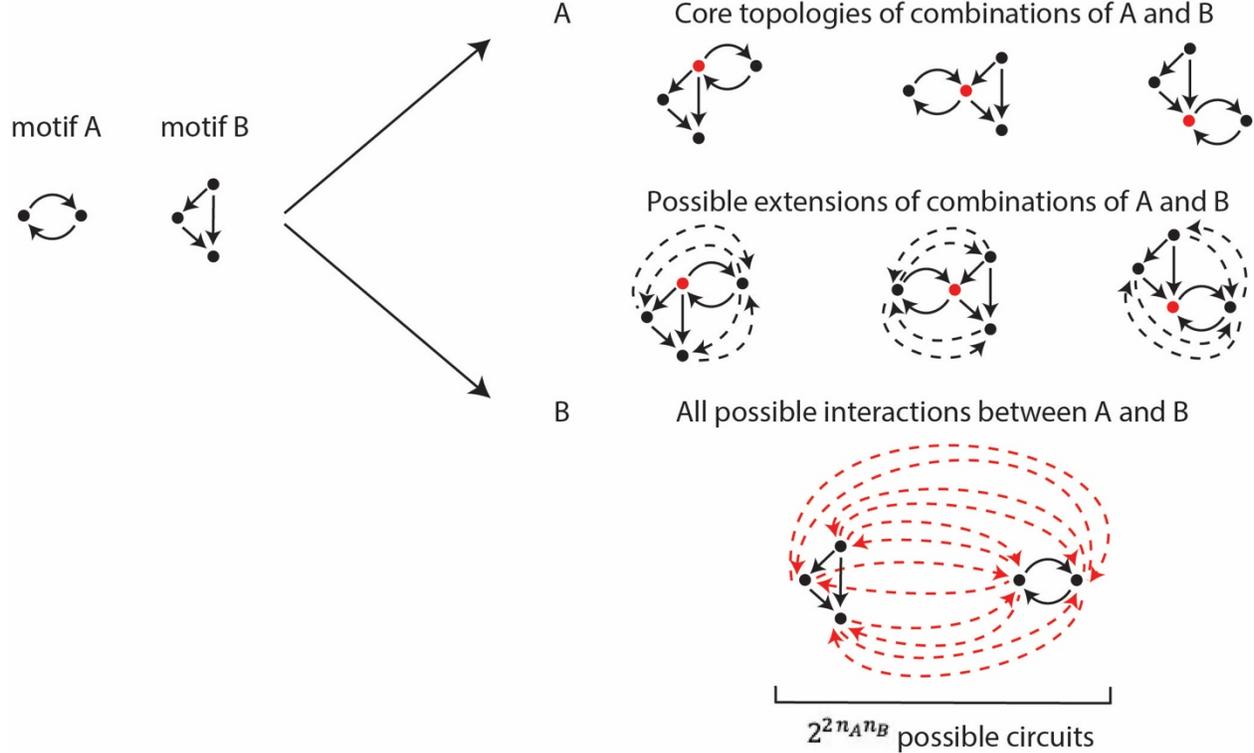

**Figure S1: A.** Definition of core topologies of combinations of motifs A and B (where they simply share at least one node), and their possible extensions in which every pair of nodes that do not participate in the same motif can be linked, exemplified for the feedback and feedforward loop motifs. The shared nodes are marked in red and the dashed edges represent the possible added links in the possible extensions. In this example there are 16 different topologies for each core combination. **B.** Example of all possible interactions between a feedback and a feedforward loop where every pair of nodes from the different motifs can interact with each other (marked in dashed red), amounting to $2^{2n_A n_B}$ possible topologies for directed networks.

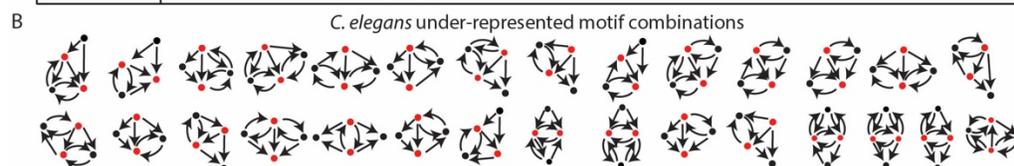

**Figure S2: A.** A table with the networks analyzed in the main text and their over-represented motif combinations. For each over-represented motif combination, we show the standard deviation of its Jaccard index (ZJ), its p-value (p) and the adjusted q-value (q). We also plot for each core topology of enriched motif combination a histogram for the frequencies of the core topology and its possible extensions in the network. Unless specified otherwise, the most frequent topology is the core topology of the motif

combinations. **B.** Topologies of additional under-represented motif combinations of the *C. elegans* neuronal network.

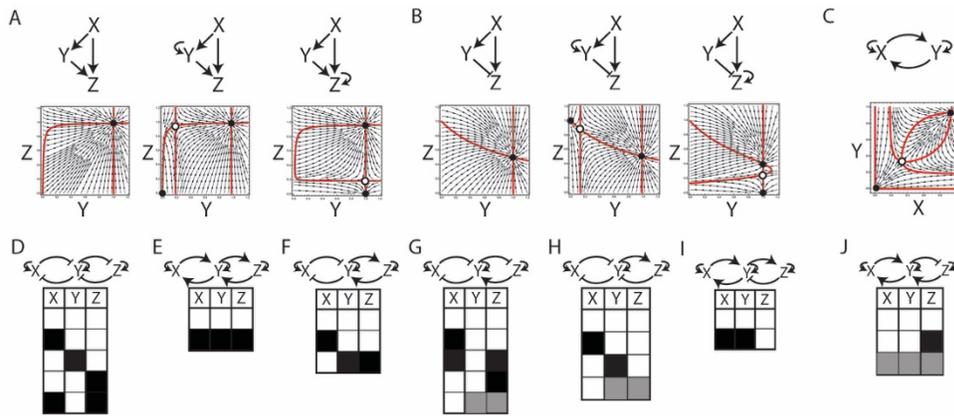

**Figure S3: A-B.** Phase portraits of the type 1 coherent (A) and incoherent (B) FFL with no autoregulation, and with a positive autoregulation of the intermediate node Y or the output node Z. **C.** Phase portrait of the lock-ON circuit in which both X and Y positively autoregulate their own levels. **D-J.** Combinations of two feedback circuits and their possible steady state where a white square represents the OFF state, black square represents an ON state, and a gray square represents an intermediate state with damped oscillations.

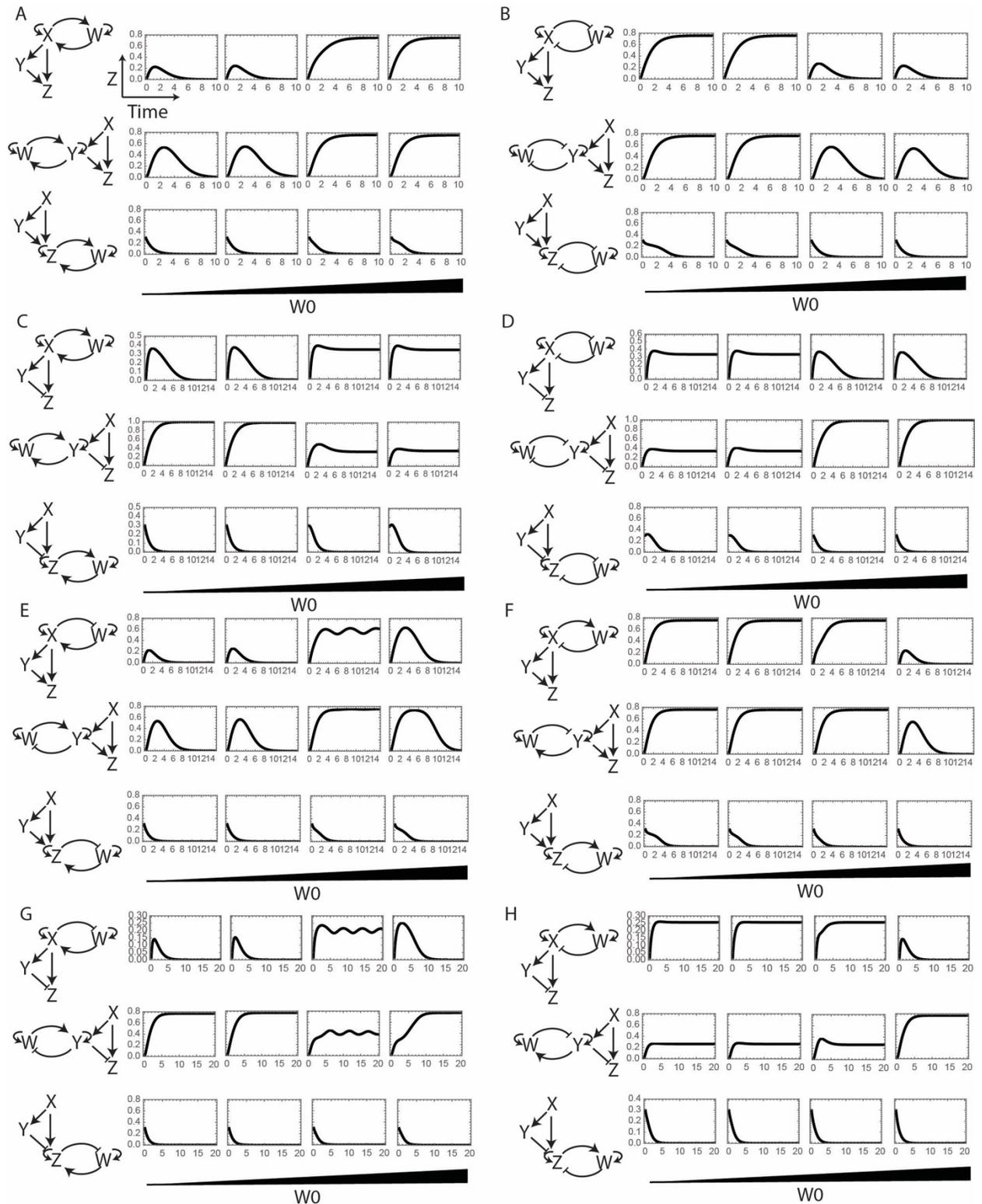

**Figure S4: A-H.** All possible core topologies of combinations of a feedforward loop and a feedback circuit where the dynamics of the FFL's output node, Z, is plotted for varying initial conditions of the node that participate only in the feedback circuit, W. Details about models and parameters are in the SI.

**Supplementary Information**

**Modeling combinations of feedforward loop circuits and autoregulation as a two-dimensional system**

To better characterize the differences between FFL circuits in which the input, intermediate or the output nodes autoregulate their own levels, we consider a model in which the input variable (X) is not a dynamical variable but rather assume that it rises in a step function manner at t=0 to its final level, $X_f$. This assumption allows us to model these circuits with only two variables, Y and Z, and thus to plot their two-dimensional phase portraits (fig. S3A-B). The model we use for combinations of a coherent FFL and positive autoregulation is therefore:

S1. $Y' = X_f - Y$

S2. $Z' = X_f Y^{n_{yz}}/(Y^{n_{yz}} + k_{yz}^{n_{yz}}) - Z$

To combine the FFL circuit with autoregulation we multiply the production term of Y or Z with $Y^{n_{yy}}/(Y^{n_{yy}} + k_{yy}^{n_{yy}})$ or $Z^{n_{zz}}/(Z^{n_{zz}} + k_{zz}^{n_{zz}})$. In Figure S3A we used $X_f = 1, k_{yz} = 0.01, k_{yy} = k_{zz} = 0.3, n_{yz} = 1, n_{yy} = n_{zz} = 3$.

The model for combinations of an incoherent FFL and positive autoregulation:

S3. $Y' = X_f - Y$

S4. $Z' = X_f Y^{n_{yz}}/(k_{yz}^{n_{yz}} + k_{yz}^{n_{yz}}) - Z$

And we combine autoregulation in the same way as described for the coherent FFL circuits. In Figure S3B we used $X_f = 1, k_{yz} = 1, k_{yy} = k_{zz} = 0.25, n_{yz} = 1, n_{yy} = n_{zz} = 3$

**Models and parameters used in Figure S3 and S4**

As described in the main text and in the Methods section, to model compositions of network motifs we consider that each node has a linear removal term. Every directed edge from X to Y is considered with the following function that multiplies the production term of Y: $X^{n_{xy}}/(k_{xy}^{n_{xy}} + k_{xy}^{n_{xy}})$ for a positive interaction and $k_{xy}^{n_{xy}}/(k_{xy}^{n_{xy}} + k_{xy}^{n_{xy}})$ for a negative interaction. We used in the Figure S3 the following parameters:

- Figure S3C-F: $k_{ij} = 0.3, n_{ij} = 3$
- Figure S3G-J: $k_{xx} = k_{ij,i \neq j} = 0.3, k_{yy} = k_{zz} = 0.2, n_{ij} = 3$
- Figure S4A-B:
  $k_{ii} = k_{wi} = k_{iw} = k_{xz} = 0.3, k_{xy} = k_{yz} = 0.01, n_{ii} = n_{wi} = n_{iw} = 3, n_{xy} = n_{xz} = n_{yz} = 1$
- Figure S4C-D:
  $k_{ii} = k_{wi} = k_{iw} = 0.3, k_{xy} = k_{xz} = 0.01, k_{yz} = 0.5, n_{ii} = n_{wi} = n_{iw} = 3, n_{xy} = n_{xz} = n_{yz} = 1$
- Figure S4E-F:
  $k_{ii} = k_{wi} = k_{xz} = 0.3, k_{iw} = 0.5, k_{xy} = k_{yz} = 0.01, n_{ii} = n_{wi} = n_{iw} = 3, n_{xy} = n_{xz} = n_{yz} = 1$
- Figure S4G-H:
  $k_{ii} = k_{wi} = k_{xz} = 0.3, k_{iw} = 0.5, k_{xy} = 0.01, k_{yz} = 0.5, n_{ii} = n_{wi} = n_{iw} = 3, n_{xy} = n_{xz} = n_{yz} = 1$